\documentclass[a4paper]{amsart}

\usepackage[utf8]{inputenc}
\usepackage{mathtools}
\usepackage{amssymb}
\usepackage{color}
\usepackage{csquotes}
\MakeOuterQuote{"}
\usepackage[sortcites=true,giveninits=true,url=false,isbn=false,doi=false]{biblatex}
\renewbibmacro{in:}{}
\AtBeginBibliography{\small}
\usepackage{hyperref}

\newtheorem{theorem}{Theorem}

\newtheorem{observation}[theorem]{Observation}
\newtheorem{lemma}[theorem]{Lemma}
\newtheorem{corollary}[theorem]{Corollary}
\newtheorem*{fact*}{Fact}

\theoremstyle{definition}
\newtheorem{definition}[theorem]{Definition}

\theoremstyle{remark}
\newtheorem{remark}[theorem]{Remark}

\newcommand{\aint}{I}

\newcommand{\Lip}{\operatorname{Lip}}

\newcommand{\dom}{\operatorname{dom}}

\newcommand{\spt}{\operatorname{spt}}
\newcommand{\dist}{\operatorname{dist}}

\renewcommand{\d}{\, \mathrm{d}}

\newcommand{\N}{\mathbb N}
\newcommand{\R}{\mathbb R}

\begin{document}
	
	\date{}
	\author{David Bate}
	\email{david.bate@warwick.ac.uk}
	\address{Zeeman Building, University of Warwick, Coventry CV4 7AL. ORCiD: \href{https://orcid.org/0000-0003-0808-2453}{0000-0003-0808-2453}}
	
	\begin{abstract}
		We show that the existing theory of Lipschitz differentiability spaces implies that a metric measure space satisfying a qualitative analogue of the WALA condition of David--Semmes has many independent Alberti representations,
		and provides conditions for when such a space is rectifiable.		
	\end{abstract}

	\title[]{Qualitative WALA, Alberti representations and Lipschitz differentiability}
	
	\keywords{}

	\thanks{This note contains no AI generated content.}
	
	\maketitle
	Recently Merlo \cite{merlo} gave a proof that the WALA condition of David--Semmes implies uniform rectifiability for measures in Euclidean space.
	The first step of the proof is to show a qualitative analogue of this statement.
	The WALA condition is a quantitative notion of differentiability and hence its qualitative analogue is very close to regular differentiability.
	In this note, we show that Merlo's qualitative statement follows by applying results from \cite{structure}, that construct independent Alberti representations in a Lipschitz differentiability space, and the results in \cite{rect,lower-den} to deduce rectifiability.

	\begin{definition}
	Fix a metric measure space $(X,d,\mu)$ throughout this note
	(a complete and separable metric space equipped with a locally finite Borel regular outer measure).
	
	Building upon the theory initiated by Cheeger \cite{cheeger}, $(X,d,\mu)$ is a \emph{Lipschitz differentiability space} if there exist countably many $U_i$, covering $X$, and Lipschitz $\varphi_i\colon X\to\R^{n_i}$ such that, every Lipschitz $f\colon X\to\R$ is differentiable at $\mu$-a.e.\ $x_0\in U$ w.r.t.\ $\varphi$: there exists a unique $Df(x_0)\in\R^n$ s.t.
		\begin{equation}
		\label{cheeger-diff}
		f(x)-f(x_0)=Df(x_0)\cdot(\varphi(x)-\varphi(x_0))+ o(d(x,x_0)).
		\end{equation}

	A \emph{curve fragment} is a bi-Lipschitz $\gamma\colon \dom\gamma \subset [0,1]\to X$ (with $\dom\gamma$ compact).
	The set of curve fragments, $\Gamma(X)$, is equipped with the Hausdorff distance between the graphs of fragments (see \cite[Definition 2.1]{structure}).
	An \emph{Alberti representation} of $\mu$ consists of a finite Borel measure on $\Gamma(X)$ for which
	\[\mu\ll\int_{\Gamma(X)}\mathcal H^1|_{\gamma}\d\eta(\gamma).\]
	(This formula always defines a Borel measure on $X$.)
	Alberti representations $\eta_1,\ldots,\eta_n$ are \emph{independent} if there exists a Lipschitz $\varphi\colon X\to\R^n$ and linearly independent cones $C_1,\ldots,C_n\subset \R^n$ s.t.\ $(\varphi\circ\gamma)'(t) \in C_i\setminus \{0\}$ for a.e.\ $t\in\dom\gamma$ and $\eta_i$-a.e.\ $\gamma\in\Gamma(X)$, for each $1\leq i \leq n$.
	This is a metric space adaptation of the ideas of Alberti \cite{alberti} and the subsequent work of Alberti, Cs\"ornyei, Preiss \cite{ACP1,ACP2},

	Following Dorronsoro \cite{dorronsoro}, David--Semmes \cite{DS} and Merlo \cite{merlo}, for $x\in \spt \mu$, $r>0$ and $f\colon X\to\R$ Borel, define $\aint(f,x,r)$ to be the average integral of $|f|/r$ over $B(x,r)$.
	If $X=\R^m$, let $\Omega(f,x,r)$ denote the infimum, over all affine $A\colon X\to\R$ of $\aint(f-A,x,r)$.
	We say $\mu$ satisfies \emph{qualitative WALA} if, for every Lipschitz $f\colon X\to\R$ and $\mu$-a.e.\ $x\in\R^m$, $\Omega(f,x,r)\to 0$ as $r\to 0$.
	\end{definition}

	In \cite{structure} several characterisations of Lipschitz differentiability spaces are given in terms of independent Alberti representations.
	The underlying idea is to first show that the existence of $n$ independent Alberti representations is characterised by a certain class $\widetilde A(\varphi)$ of subsets of $X$ being $\mu$-null.
	Secondly, show that any $\widetilde A(\varphi)$ set is $\mu$-null by constructing a Lipschitz function differentiable almost nowhere on such a set.
	In this note we show that the same results may be applied to measures satisfying qualitative WALA.

	In a Lipschitz differentiability space, one often considers the equivalent (up to a countable decomposition) condition that there is a bound on the dimension of any ``infinitesimally linearly independent'' vector valued Lipschitz function:
	If $F\colon X\to\R^{d}$ is Lipschitz and, for a set of $x_0$ with positive measure, $\Lip(a(x_0)\cdot F,x_0)>0$ for every $a\in\mathbb S^{d-1}$, then $d\leq n$.
	This idea originated in \cite{cheeger} and has been reused extensively, see e.g.\ \cite{keith,primer,kakeya}.
	A consequence is that no regularity of $\varphi$ is necessary: if $\varphi$ is any \emph{function} satisfying \eqref{cheeger-diff} for all $f$, then there is another Lipschitz $\varphi$ also satisfying \eqref{cheeger-diff} for all $f$.
	The same idea can be applied to qualitative WALA.

	\begin{observation}
		\label{indep}
		For $d\geq n+2$, let $f_1,\ldots,f_{d} \colon X\to\R$, $\varphi\colon X \to \R^n$ Borel, $x\in X$ and $r>0$.
		For $\epsilon>0$ and each $1\leq i \leq d$ suppose that there exist $A_i\in\R^n$ such that
		\(\aint(f_i-A_i\cdot \varphi,x_0,r) \leq \epsilon.\)
		Write $F\colon X\to\R^d$ for the function whose components are the $f_i$.
		Then there exists $a\in\mathbb S^{d-1}$ such that
		\(\aint(a\cdot (F-F(x_0)),x_0,r) \leq d\epsilon.\)
		This motivates the following definition.
		For $\mathcal F$ a set of Borel $g\colon X\to\R$, define $\Omega_\mathcal F(f,x,r)$ similarly to $\Omega$, but using elements of $\mathcal F$ rather than affine functions.

		We combine this observation with a doubling condition:
		suppose $D\subset X$ and $C,T>0$ such that
		\begin{equation}
			\label{doubling}
			\mu(B(x,2r)) \leq C \mu(B(x,r)) \quad \forall x\in D, 0<r<T.
		\end{equation}
		For $r<T/4$ set $g=a_r\cdot(F-F(x_0))/r$.
		If $L>0$, $F$ is $L$-Lipschitz and $g(x)\geq \epsilon r$, then $g$ is bounded below by $\epsilon r/2$ on $B(x,\eta r)$ for $\eta=\min(1,\epsilon/2L)$.
		If $x,x_0\in D$ then $\mu(B(x,\eta r))\gtrsim_{L,\eta} \mu(B(x_0,2r))$ and so $\aint(a_r\cdot (F-F(x_0)),x_0,r)\to 0$ implies
		\begin{equation}
			\label{wala-diff}
			\lim_{r\to 0}\sup_{B(x_0,r)\cap D}|a_r\cdot(F-F(x_0))|/r \to 0.
		\end{equation}
	\end{observation}
	Thus, the only significant difference between qualitative WALA and differentiability is that the linear combination in \eqref{wala-diff} may depend on $r$, whereas for differentiability it does not.
	In \cite{BS,structure} Lipschitz functions are constructed with both large and arbitrarily small difference quotients, which implies that they are not differentiable.
	The uniqueness of the linear combination does not play a role in those arguments and so the methods are easily adapted to the qualitative WALA setting.
	
	We require modifications of two statements from \cite{structure} that record more properties of each proof. 
	\begin{lemma}[\cite{structure} Lemma 4.2]
		\label{lemr1}
	Let $(X,d,\mu)$ be a metric measure space, $S\subset X$ Borel and $0<\alpha<1$.
	Suppose that there exist a sequence of 1-Lipschitz functions $f_i\colon X\to \R$, and sequences $\theta_i>0$ and $\sigma_i>0$ such that, for every $x_0\in S$ and any $y,z\in B(x_0,\sigma_i)$,
			\[|f_i(y)-f_i(z)| \leq \theta_i d(y,z).\]
	Suppose also that there exist decreasing sequences $R_i>r_i>0$ such that $R_{i+1}<\sigma_i$, $R_{i+1}<2^{-(i+1)}\alpha r_i$, $\theta_i\leq \alpha 2^{-i}$ and $|f_i|_\infty\leq R_i$ for every $i\in\N$.
	Then there exists a 3-Lipschitz function $F\colon X\to \R$ such that
		for almost every $x_0\in S$, each $i\in\N$ and any $x\in X$ with $r_{i}\leq d(x,x_0)\leq R_{i}$,
			\begin{equation}
				\label{fun2}
				\frac{|F(x)-F(x_0)-(f_{i}(x)-f_{i}(x_0))|}{d(x,x_0)} \leq \alpha.
			\end{equation}
	\end{lemma}
 
	\begin{lemma}[\cite{structure} Lemma 6.5]
		\label{lemr2}
		Let $(X,d,\mu)$ be a metric measure space, $\varphi \colon X\to\R^n$ 1-Lipschitz, $S\in \widetilde A(\varphi)$ and $U\supset S$ compact.
		Then for any $0<\delta,\xi,\lambda<1$ there exists a $R>0$ such that, for any $0<r<R$, there exists a compact $S'\subset S$ with $\mu(S')\geq \mu(S)-\xi$, $\sigma>0$ and a 1-Lipschitz function $f\colon X\to\R$ such that
		\begin{enumerate}
			\item For each $x_0\in S'$ there exists $v\in\mathbb S^{n-1}$ such that, for every $x\in U$ with $r \leq d(x,x_0)\leq R$ and $|(\varphi(x)-\varphi(x_0)) \cdot v|\geq \lambda d(x,x_0)$,
				\[|f(x)-f(x_0)|\geq |(\varphi(x)-\varphi(x_0)) \cdot v|/2.\]
			\item For every $x_0\in S'$ and $y,z\in B(x_0,\sigma)$
				\[|f(y)-f(z)|\leq \delta d(y,z).\]
		\end{enumerate}
	\end{lemma}

	In \cite[Lemma 2.1]{BS} a condition similar to \eqref{BS1} below is shown to be equivalent to the uniqueness of the Cheeger derivative at $x_0\in X$ with respect to $\varphi\colon X \to \R^{n}$.
	By finding points given by this condition in many disjoint annuli, the proof of \cite[Theorem 2.4]{BS} constructs many infinitesimally linearly independent Lipschitz functions.
	The main contribution of this note is that, if these points appear at the same scale, then the functions constructed violate Observation \ref{indep}.

	\begin{definition}[Compare to \cite{structure} Definition 3.6]
		\label{structured}
		A compact set $U\subset X$ has \emph{structured neighbourhoods} if the following are true. There exists a compact $D\supset U$, $T,C>0$ such that \eqref{doubling} holds.
		For $d\in\N$ there exist $M_1,\ldots,M_d\in\N$ and $N_1,N_2,N_3\in \N$ with $M_j\geq M_{j-1}+N_3+ N_1+1$ and $0<\lambda,\rho<1$, an increasing sequence $i(k)\to\infty$ and a finite set $\mathcal P$ of 1-Lipschitz functions $\varphi\colon X\to\R^n$ such that the following holds:
		For each $x_0\in U$ and $k\in\N$, $1\leq j\leq d$ there exists $\varphi\in\mathcal P$ and $m_1,m_2\in\N$ with $N_3\geq m_2\geq m_1+1+N_2$ such that, for each $l=1,2$,
		\begin{equation}
		\label{BS1}
		\sup \frac{|(\varphi(x)-\varphi(x_0))\cdot v|}{d(x,x_0)} \geq \lambda \|v\| \quad  \forall v\in \R^n,
	\end{equation}
	where the supremum is over ${x\in D\cap B(x_0,\rho^{i(k)+M_j+m_l})\setminus B(x_0,\rho^{i(k)+M_j+m_l+1})}$.
	\end{definition}

	\begin{observation}
		\label{annuli}
		Let $(X,d,\mu)$ be a metric measure space and suppose that $U\subset X$ has structured neighbourhoods and $S\in\widetilde A(\varphi)\cap U$ for each $\varphi\in \mathcal P$.
	 	Suppose those parameters satisfy $\rho^{N_2}<\lambda/8$, $\rho^{N_1}<(\rho\lambda/48d)^d/d$.
		Then for any $\delta>0$ there exists a finite collection $\mathcal S$ of Lipschitz functions $\Phi\colon X\to\R^d$, a $S'\subset S$ with $\mu(S')\geq \mu(S)-\delta$ and a (un-relabeled) subsequence of $i(k)$ such that, for each $x_0\in S'$ and $k\in\N$,
		\begin{equation}\label{max}
			\max_{\Phi\in\mathcal S}\min_{a\in\mathbb S^{d-1}} \sup_{D\cap B(x_0,\rho^{i(k)})}|a\cdot (F_\Phi-F_\Phi(x_0))| \rho^{-i(k)} \gtrsim_{C,\rho,M_d,N_1,N_3,d} 1.
		\end{equation}
	\end{observation}

	\begin{proof}
		Set $\theta_k\leq \alpha 2^{-k}$ for each $k\in\N$.
		For each $\varphi\in\mathcal P$ and $k\in\N$, set $\delta=\theta_{k}$, $\xi=2^{-k}\mu(S)/2|\mathcal P|$ and let $S'_{\varphi,k}$, $R_{\varphi,k}$ and $v_{\varphi,k}$ be obtained from Lemma \ref{lemr2} applied to $\varphi$.
		By taking a subsequence, we may suppose that $\rho^{i(k)}<\min_\varphi,R_{\varphi,k}$ for each $k\in \N$.
		Now set $R_k=\rho^{i(k)}$ and $r_k=\rho^{i(k)+M_d+N}$ for each $i\in\N$ and $0<\alpha<\min(\lambda\rho/16,\rho^{N_1})$.
		Obtain $f_{\varphi,k}$ and $\sigma_{\varphi,k}$ from Lemma \ref{lemr2}.
		Then $f_{\varphi,k}$ has the property that, for every $x_0\in S'_\varphi:=\cap_{k}S'_{\varphi,k}$ satisfying Definition \ref{structured} for $i(k)$, $\varphi$, and $M_j$, there exists $m_1,m_2\in\N$ with $N_3\geq m_2\geq m_1+1+N_2$ and $x_1,x_2\in D$ with $\rho^{i(k)+M_j+m_l+1}\leq d(x_l,x_0) \leq \rho^{i(k)+M_j+m_l}$ and
		\begin{equation}
		\label{big-fj}\frac{|f_{\varphi,k}(x_l)-f_{\varphi,k}(x_0)|}{d(x,x_0)} \geq \frac{\lambda}{2}.
		\end{equation}

		We now modify each $f_{\varphi,k}$ similarly to the proof of \cite[Lemma 4.1]{structure}.
		For each $1\leq j\leq d$, post-compose $f_{\varphi,k}$ with the function $t\mapsto \dist(t,2\rho^{{i(k)+M_j}}\mathbb Z)$ to obtain $g_{\varphi,k}^j$.
		This composition can break \eqref{big-fj}, for possibly one of $x_1$ or $x_2$, but not both.
		Indeed, if $|g_{\varphi,k}^j(x_2)-g_{\varphi,k}^j(x_0)|\leq \lambda d(x_2,x_0)/4$ then $\dist(f_{\varphi,k}(x_0),\mathbb Z\rho^{i(k)+M_j})\leq \rho^{i(k)+M_j+m_2}$ and
		\begin{align*}|g_{\varphi,k}^j(x_1)-g_{\varphi,k}^j(x_0)|&\geq \lambda\rho^{i(k)+M_j+m_1+1}/2-2\rho^{i(k)+M_j+m_2}\\&\geq \rho^{i(k)+M_j+m_1+1}(\lambda/2-2\rho^{N_2})\\&\geq \lambda \rho d(x_1,x_0)/8.\end{align*}
		We write $x'$ for the point $x_1,x_2$ satisfying this inequality.
		On the other hand, for any $d(x,x_0)\geq \rho^{i(k)+M_{j-1}+N_3}$,
		\begin{equation*}
		\frac{|g_{\varphi,k}^j(x)-g_{\varphi,k}^j(x_0)|}{d(x,x_0)} \leq \frac{\rho^{i(k)+M_j}}{\rho^{i(k)+M_{j-1}+1}}=\rho^{N_1}.
		\end{equation*}
		By taking a subsequence of $i(k)$ if necessary, we may suppose that $i(k+1)\leq \min_{\varphi}\rho_{\varphi,k}$ for each $k\in\N$ and that the hypotheses of Lemma \ref{lemr1} are satisfied (for the sequence $i(k)$).
		Finally, set $S'=\cap_{\varphi,k}S'_{\varphi,k}$, so that $\mu(S')\geq \mu(S)/2$.

		Now fix $\Phi:=(\varphi_1,\ldots,\varphi_d)\in \mathcal P^d$.
		For each $1\leq j\leq d$ apply Lemma \ref{lemr1} to the sequence $g^j_{\varphi_j,k}$ to obtain a Lipschitz function $F_\Phi^j$.
		Write $F_\Phi$ for the function whose coordinates are $F_\Phi^j$.
		For some $k\in\N$ suppose that $x_0\in S'$, $\varphi_j$ and $M_j$ satisfy Definition \ref{structured} for each $1\leq j\leq d$.
		Then by the above, for each $1\leq j\leq d$ and any $x\in D$ with $d(x,x_0)\geq \rho^{i(k)+M_{j-1}+N_3}$,
		\begin{equation}
			\label{smallF}
		\frac{|F_{\Phi}^j(x)-F_{\Phi}^j(x_0)|}{d(x,x_0)} \leq \rho^{N_1} + \alpha \leq 2\rho^{N_1},
		\end{equation}
		and there exist $x_j\in D$ with $\rho^{i(k)+M_j+N_3}\leq d(x_j,x_0)\leq \rho^{i(k)+M_j}$ such that
		\begin{equation}
			\label{largeF}
			\frac{|F_\Phi^j(x'_j)-F_\Phi^j(x_0)|}{d(x'_j,x_0)} \geq \frac{\lambda\rho}{8}-\alpha\geq \frac{\lambda\rho}{16}.
		\end{equation}
		For any $a\in\mathbb R^{d}$, there exists $z\in D\cap B(x_0,\rho^{i(k)})\setminus B(x_0,\rho^{i(k)+M_d+N_3})$ with
		\begin{equation}
			\label{combination}
			|a\cdot(F_{\Phi}(z)-F_{\Phi}(x_0))|\geq \|a\|\rho^{N_1} d(z,x_0).
		\end{equation}
		Indeed, suppose this is not true.
		For any $1\leq i\leq d$, applying \eqref{largeF} for $i$, \eqref{smallF} for $j>i$, the fact that each $F_\Phi^j$ is 1-Lipschitz for $j<i$, the negation of \eqref{combination} and the triangle inequality imply
		\[\lambda\rho|a_i|/16 < d \max_{j<i}|a_j| + 2d\rho^{N_1}\max_{j>i}|a_j| + d\rho^{N_1} \max_{j}|a_j|\leq d \max_{j<i}|a_j| + 3 d\rho^{N_1}\max_{j}|a_j|.\]
		Iterating over $j$ from $1$ to $d$ gives $\max |a_j| < d(48d/\rho\lambda)^d\rho^{N_1}\max|a_j|$, a contradiction.
	\end{proof}

	\begin{theorem}
		\label{theo}
		Let $(X,d,\mu)$ be a metric measure space and suppose that $U\subset X$ has structured neighbourhoods with parameters satisfying the conditions in Observation \ref{annuli}.
		Further, for any Lipschitz $F\colon X \to \R^{d}$, every $r>0$ and almost every $x_0\in U$, there exists $a_r\in\mathbb S^{d-1}$ such that \[\lim_{r\to 0} \sup_{D\cap B(x_0,r)} |a_r\cdot (F-F(x_0))|r^{-1}= 0.\]
	Then $\mu|_U$ has a countable decomposition into pieces with $n$-independent Alberti representations (each with respect to some $\varphi\in\mathcal P$).
	\end{theorem}

	\begin{proof}
		The proof closely follows the proof of \cite[Theorem 6.6]{structure}.
		By \cite[Proposition 2.10]{perturbations}, there exists a countable Borel decomposition $U=S\cup \cup_j U_j$ where each $\mu|_{U_i}$ has $n$-independent Alberti representations with respect to some $\varphi\in\mathcal P$, and $S\in\widetilde A(\varphi)$ for all $\varphi\in\mathcal P$.
		Apply Observation \ref{annuli} for $\delta=\mu(S)/2$ to obtain $S'$ and $\mathcal S$.
		Since $\mathcal S$ is finite, for any $x_0\in S$ a single $\Phi\in\mathcal S$ attains the maximum in \eqref{max} for infinitely many $k$, contradicting the hypotheses.
		Consequently, $\mu(S')=\mu(S)=0$.		
	\end{proof}

	\begin{remark}
		The hypotheses of Definition \ref{structured} can be weakened and still arrive to Theorem \ref{theo}; we chose the hypotheses as to require the least modification of the proof given in \cite{structure}.
	\end{remark}

	The article \cite{rect} characterises when a measure with many independent Alberti representations, in particular a chart $U$ in a Lipschitz differentiability space, is rectifiable.
	The sufficient conditions were improved in \cite{lower-den}.
	Applying these results gives the following.
	\begin{corollary}
		In addition to the hypotheses in Theorem \ref{theo}, suppose that
		\[\Theta^{*,n}(\mu,x):= \limsup_{r\to 0} \frac{\mu(B(x,r))}{r^n} >0\]
		for $\mu$-a.e.\ $x\in U$.
		Then $\mu|_U$ is $n$-rectifiable.
	\end{corollary}

	The conditions in Definition \ref{structured} can be explicitly verified in Euclidean space.

	\begin{corollary}[\cite{merlo} Theorem 1.1]
		For $\alpha\geq 0$, let $\mu$ be a Borel measure on $\R^m$ such that $\Theta^{*,\alpha}(\mu,x)<\infty$ and
		\[\Theta_{*}^\alpha(\mu,x):=\liminf_{r\to 0} \frac{\mu(B(x,r))}{r^\alpha} >0\]
		for $\mu$-a.e.\ $x\in \R^m$.
		Suppose that $\mathcal F$ is a finite dimensional subspace of Borel functions $f\colon \R^m\to\R$ and that, for every Lipschitz $f\colon \R^m\to \R$, $\Omega_{\mathcal F}(f,r,x)\to 0$ as $r\to0$ for $\mu$-a.e.\ $x\in \R^m$.
		Then $\alpha=n\in\N$ and $\mu$ is $n$-rectifiable.
	\end{corollary}

	\begin{proof}
		Let $n=\lceil\alpha\rceil$ and $d\geq \dim\mathcal F+2$.
		We must check that the hypotheses of Definition \ref{structured} and Theorem \ref{theo} hold.
		There exists a countable decomposition $\R^m=N\cup \bigcup_i D_i$ with each $D_i$ compact, $\mu(N)=0$ and, for each $i\in \N$, $x\in D_i$,
		\[\frac{r^\alpha}{i} \leq \mu(B(x,r)) \leq i r^\alpha \quad \forall 0<r< \frac{1}{i}.\]
		Thus there exists $\lambda,\rho>0$ such that, for each $x_0\in D_i$ and $k\in\N$ with $\rho^k\leq 1/i$, \eqref{BS1} holds for $i(k)=k, m_l=0$ and $\varphi$ an orthogonal projection (an Ahlfors $\alpha$-regular measure cannot be close to an $n-1$ dimensional subspace, quantitatively).
		By reducing $\lambda>0$ slightly, these $\varphi$ may be chosen from a finite prescribed set $\mathcal P$.
		For any $N_2\in\N$, provided $N_3$ is chosen large enough, the compactness of $\mathbb S^{m-1}$ gives the existence of $m_1,m_2$ satisfying \eqref{BS1}.
		Therefore, Theorem \ref{theo} gives a countable Borel decomposition of each $D_i$ into sets on which $\mu$ has $n$-independent Alberti representations.

		The number of independent Alberti representations a measure can possess is bounded above by its Hausdorff dimension.
		This is a consequence of Guth's multilinear Kakeya inequality \cite{guth} or the structural theorem of De Philippis--Rindler \cite{DPR}.
		See \cite[Theorem 5.3]{kakeya} or \cite[Theorem 2.17]{perturbations} for details respectively.
		Thus we must have $\alpha=n$, and moreover, $\mu$ is $n$-rectifiable.
	\end{proof}

	\begin{remark}
	If $X=\R^m$ for some $m\in\N$, one may instead use \cite[Theorem 2.18]{perturbations} to deduce rectifiability (as is done in \cite[Proposition 3.1]{merlo}).
	The proof combines the Besicovitch--Federer projection theorem and De Philippis--Rindler \cite{DPR}.
	\end{remark}

	\begin{remark}
		When $d=n+1$ and porous sets have measure zero, the Alberti representations obtained from Theorem \ref{theo} are universal.
		Hence, if $X$ is covered by countably many such $U$, $X$ is a Lipschitz differentiability space by \cite[Theorem 7.8]{structure}.
		Conversely, any $n$-dimensional chart in a Lipschitz differentiability spaces satisfies the hypotheses of Definition \ref{structured}.
		This follows from the existence of $n$-independent Alberti representations \cite[Theorem 6.6]{structure} and that porous sets have measure zero \cite[Theorem 2.4]{BS}.
	\end{remark}

	\printbibliography
\end{document}